\documentstyle[12pt]{article}
\newcommand{\ptl}{\partial}

\newcommand{\G}{\Gamma}

\newcommand{\vs}{\varsigma}

\newcommand{\ol}{\overline}

\newcommand{\kn}{\mbox{ker}}

\def\sF{\hbox{$\sc I\hskip -2.5pt F\!$}}

\def \Z{\hbox{$Z\hskip -5.2pt Z$}}

\def \C{\hbox{$C\hskip -5pt \vrule height 6pt depth 0pt \hskip 6pt$}}

\def\qed{\hfill \hfill \ifhmode\unskip\nobreak\fi\ifmmode\ifinner
         \else\hskip5pt\fi\fi
 \hbox{\hskip5pt\vrule width4pt height6pt depth1.5pt\hskip 1 pt}}
\def\a{\alpha}
\def\b{\beta}
\def\d{\delta}
\def\D{\Delta}
\def\g{\gamma}
\def\G{\Gamma}
\def\l{\lambda}

\def\sc{\scriptstyle}
\def\ssc{\scriptscriptstyle}
\def\F{\hbox{$I\hskip -4pt F$}}\def \Z{\hbox{$Z\hskip -5.2pt Z$}}
\def\dis{\displaystyle}
\def\cl{\centerline}

\def\ol{\overline}

\def\rar{\rightarrow}
\def\Rar{\Rightarrow}

\def\Lar{\Leftarrow}

\def\Lra{\Leftrightarrow}
\def \II{\hbox{$\cal I$}}
\def \KK{\hbox{$\cal K$}}
\def \LL{\hbox{$\cal L$}}
\def \MM{\hbox{$\cal M$}}
\def\bs{\backslash} \def\vsp{{}}

\def\vs{\vspace*}

\begin{document}
\par\
\par\
\par
\cl
{{\bf SIMPLE ALGEBRAS OF WEYL TYPE$\sc\,$}\footnote
{AMS Subject Classification - Primary:  17B20, 17B65, 17B67, 17B68.\\
\indent \hskip .3cm  This work is supported by NSF of China and a Fund from
National Education Department of China.}}
\vsp
\par\
\vs{-7pt}
\par
\centerline{Yucai Su$^\dag$$\ssc\,$\footnote
{This author was partially supported by  Academy of Mathematics and Systems Sciences
during his visit to this academy.}
 \,\,\,\, and \,\,\,\, Kaiming Zhao$^\ddag$}
\vsp
\par\
\vs{-5pt}
\par
 $^\dag$Department of Applied Mathematics,  Shanghai Jiaotong University,
Shanghai 200030, P. R. China

$^\ddag$Institute of Mathematics,  Academy of Mathematics and systems Sciences,
Chinese \break Academy of Sciences,  Beijing 100080, P. R. China
\vsp
\par\
\par

{\bf Abstract} Over a field $\F$ of any  characteristic,
for a commutative associative algebra
$A$ with an identity element and for the polynomial algebra $\F[D]$ of
a commutative derivation subalgebra
$D$ of $A$,
the associative and the Lie algebras
of Weyl type on  the same vector space $A[D]=A\otimes\F[D]$
are defined. It is proved that $A[D]$, as a Lie algebra
(modular its center) or as an associative
algebra, is simple if and only if
$A$ is $D$-simple and $A[D]$ acts faithfully on $A$. Thus
a lot of  simple algebras are obtained.
\par\
\par
For a  long time, simple Lie algebras and simple associative algebras  have
been two central objects in the theory of algebra. The four well-known series
of infinite dimensional simple Lie algebras of Cartan type have played
important roles in the structure theory of Lie algebras. Generalizations of the
simple Lie algebras of Cartan type over a field of characteristic zero have
been obtained by Kawamoto [1], Osborn [2], Dokovic and Zhao [3,4,5], Osborn and
Zhao [6,7] and Zhao [8]. Passman [9], Jordan [10] studied the Lie algebras
$AD=A\otimes D$ of generalized Witt type constructed from a commutative
associative algebra $A$ with an identity element and its commutative
derivation subalgebra $D$ over a field $\F$ of arbitrary  characteristic.
Passman proved that $AD$ is simple if and only if $A$ is $D$-simple and
$AD$ acts faithfully on $A$. Xu [11] studied some of these simple Lie
algebras of Witt type and other Cartan types Lie algebras, based on the pairs
of the tensor algebra of the group algebra of an additive subgroup of ${\F}^n$
with the polynomial algebra in several variables and the subalgebra of
commuting locally finite derivations. Su, Xu and Zhang [12] gave the structure
spaces of the generalized simple Lie algebras of Witt type constructed in [11].
We\footnote{ Su Y, Zhao K. Second cohomology group of generalized Witt type
Lie algebras and certain representations. Submitted for publication.}
determined the second cohomology group and gave some representations of the Lie
algebras of generalized Witt type which are some Lie algebras defined by
Passman, more general than those defined by Dokovic and Zhao, and slightly more
general than those defined by Xu.

\par
Throughout this paper,  $\F$ is a  field of any characteristic.
\par
In this paper, for a commutative associative algebra
$A$ with an identity element over $\F$
and for the polynomial algebra of
its commutative derivation subalgebra $D$, we construct the associative algebras $A\otimes\F[D]$
of Weyl type and determine
 the necessary and sufficient conditions for $A\otimes\F[D]$ to be simple, as
a Lie algebra and as associative algebra  respectively, see Theorems 1.1, 1.2.  Using this construction, by giving specific $A$ and $D$, we
then obtain a large class of simple  algebras.
Another interesting observation is that
all Lie algebras mentioned above are Lie subalgebras of the Weyl type Lie
algebras constructed in this paper.
\par\
\vs{7pt}\par
{\bf 1 Simple algebras of Weyl type}
\vs{10pt}\par
Let $A$ be an  commutative associative algebra with an identity element 1
over $\F$, let $D$ be a nonzero $\F$-vector
space spanned by some commuting $\F$-derivations of $A$.
Let
$$
\{\ptl_i\,|\,i\in I\},
\mbox{ where $I$ is some indexing set},
\eqno(1)$$
be an $\F$-basis of $D$. Choose a total ordering on $I$. Denote
$$
J=\{\a{\sc\!}={\sc\!}(\a_i\,|\,i\in I)\,|\,\a_i{\sc\!}\in{\sc\!}\Z_+,
\forall\,i{\sc\!}\in{\sc\!} I
\mbox{ and  $\a_i{\sc\!}={\sc\!}0$ for all but a finite number of
$i{\sc\!}\in{\sc\!} I\}$},
\eqno(2)$$
where $\a=(\a_i\,|\,i\in I)$, we shall simply write $\a=(\a_i)$,  is a
collection of nonnegative integers
indexed by $I$.
Define a total ordering on $J$ by
$$
\a<\b\Lra |\a|<|\b|\mbox{ or }|\a|=|\b|
\mbox{ but there }\,\,\, \exists\,i\in I,\a_i<\b_i
\mbox{ and }\a_j=\b_j,\,\forall\,j\in I,j<i,
\eqno(3)$$
where $|\a|=\sum_{i\in I}\a_i$ is called the {\it level} of $\a$.
Denote by $\F[D]$ the polynomial algebra (the group algebra) of $D$
with basis
$$
B=\{\ptl^{(\a)}=\prod_{i\in I}\ptl_i^{\a_i}\,|\,\a\in J\}.
\eqno(4)$$
Let $A[D]=A\otimes\F[D]$ be the tensor product of $A$ and $\F[D]$,
which acts naturally on $A$ by
$$
u\otimes\ptl^{(\a)}:x\mapsto
u\cdot\ptl^{(\a)}(x)=
u\cdot(\prod_{i\in I}\ptl_i^{\a_i})(x),
\ \forall\,u,x\in A,\,\ptl^{(\a)}\in B,
\eqno(5)$$
where $( \ptl_1\ptl_2 \cdots \ptl_n)(x)=( \ptl_1(\ptl_2 \cdots (\ptl_n(x) )\cdots))$.
 This gives rise to a linear transformation
$$
\theta:A  [D] \rar{\rm Hom}_{\sF}(A,A).
\eqno(6)$$
For any $\a\in J$, set
$$
\matrix{
{\rm supp}(\a)=\{i\in I\,|\,\a_i\ne0\},
\vs{4pt}\hfill\cr
J(\a)=\{\g\in J\,|\,\g_i\le\a_i,\,\forall\,i\in I\},
\vs{4pt}\hfill\cr
\dis({}^\a_\g)=\prod_{i\in I}({}^{\a_i}_{\g_i}),\ \g\in J(\a).
\hfill\cr}
\eqno(7)$$
To obtain an associative algebra structure on $A[D]$ so that
$\theta$ is a homomorphism of  associative algebras, we define the product
as follows
$$
(u\otimes\ptl^{(\a)})\cdot(v\otimes\ptl^{(\b)})=
u\sum_{\g\in J(\a)}({}^\a_\g)\ptl^{(\g)}(v)\otimes\ptl^{(\a+\b-\g)},
\eqno(8)$$
for all $u,v\in A,\,\a,\b\in J$.
Then (6) defines $A$ as an $A[D]$-module.
We call the associative algebra $(A[D],\cdot)$
an {\it associative algebra of Weyl type}.
For  Weyl algebras we refer the reader to [13].
We shall simply denote $u\otimes d$ by $ud$ for $u\in A,d\in\F[D]$.
For any
$$
x=\sum_{\a\in J}u_\a\ptl^{(\a)}\in A[D],
\eqno(9)$$
we say $x$ has
{\it leading term} \,$ld(x)=u_\b\ptl^{(\b)}$, {\it leading degree}
\,$deg(x)=\b$ and {\it leading level} \,$lev(x)=|\b|$ if
$u_\b\ne0$ and for all $\a\in J,u_\a\ne0\Rar\a\le\b$. Define
$lev(0)=-\infty.$ We define the {\it support} of $x$ to be the
set $\{\a\in J\,|\,u_\a\ne0\}$.
Set
$$
\F_1=\{u\in A\,|\,D(u)=0\}.
\eqno(10)$$
We use this notation because it is a field extension of $\F$ when $A$
is {\it $D$-simple}, i.e., $A$ has no nontrivial $D$-stable ideals,
see [9].
\par
Define the binary operation $[\cdot,\cdot]$  to be the
usual induced Lie bracket on the associative algebra $(A[D],\cdot)$ so that
$(A[D],[\cdot,\cdot])$ is a Lie algebra. Obviously, $\F_1$ is contained
in the center of $(A[D],[\cdot,\cdot])$. Denote
$$
\ol{A}[D]=A[D]/\F_1,
\eqno(11)$$
whose induced Lie bracket is also denoted by $[\cdot,\cdot]$. We call the
Lie algebra $\ol{A}[D]$ a {\it Lie algebra of Weyl type}.
\par\
\par
{\bf Theorem 1.1} The Lie algebra $\ol{A}[D]$ is simple if and only if
$A$ is $D$-simple and $\F_1[D]$ acts faithfully on $A$.
\par
{\bf Proof.}
``$\Rar$'': If $\II$ is a $D$-stable ideal of $A$, then
clearly $\II[D]$ is a Lie ideal of $A[D]$.
If $\II\ne0$, then $\II[D]\not\subset\F_1$, and so
$\II[D]=A[D]$.
Thus $\II$ must be $A$.
Furthermore, the kernel $\kn\theta$ of the Lie
homomorphism $\theta$ is a Lie ideal of $A[D]$.
If $\kn\theta=A[D]$, then in particular,
$D$ acts trivially
on $A$, contradicting its definition as a nonzero subspace of
${\rm Der}_{\sF}(A)$. Thus $\kn\theta\subset\F_1$.
If $\theta(f)=0$ for $f\in A$, then $0=\theta(f)(1)=f$.  Hence,
$\kn\theta=0$, i.e., $A[D]$ acts faithfully on $A$. In
particular, $\F_1[D]$ acts faithfully on $A$.
\par
``$\Lar$'': We prove the sufficient conditions by several claims.
\par
{\bf Claim 1}. $A[D]$ acts faithfully on $A$.
\par
Suppose $\kn\theta\ne0$. Let $m\ge1$ be the minimal support
size of  nonzero elements in $\kn\theta$, i.e., there exist
\F-linearly independent
elements $d_1,\cdots,d_m$ of $B$ with $\kn\theta\cap\sum_{p=1}^mA d_p\ne0$ and
such that all nonzero elements in this intersection have all their
$A$-coefficients being nonzero. Set
$$\KK={\rm span}\{a_1\in A\,|\,v=\sum_{p=1}^m a_pd_p\in\kn\theta
\mbox{ for some }a_p\in A\}.
\eqno(12)$$
Then $\KK$ is a nonzero subspace of $A$. For any $\ptl\in D$, using (8)
and the fact that $\kn\theta$ is a Lie ideal of $A[D]$,
we have
$$\sum_{p=1}^m\ptl(a_p)d_p=[\ptl,v]\in\kn\theta,
\eqno(13)$$
where $v$ is written as in (12).
Thus by definition (12),
$\ptl(a_1)\in\KK$, i.e., $\KK$ is a $D$-stable subspace.
Obviously, $A\cdot\kn\theta\subset\kn\theta$, and so $A a_1\subset\KK$, i.e.,
$\KK$ is a $D$-stable ideal of $(A,\cdot)$. Thus $\KK=A$.
In particular $1\in\KK$ and so there exists
some $w=d_1+\sum_{p=2}^m a_pd_p\in\kn\theta$.
Then for all $\ptl\in D$, we have $\sum_{p=2}^m\ptl(a_p)d_p
=[\ptl,w]\in\kn\theta$ with support size $\le m-1$.
This means that $\ptl(a_p)=0$ for all $p=2,\cdots,m$ and all $\ptl\in D$.
In other word, $a_p\in\F_1$ and so $0\ne w\in\F_1[D]\cap\kn\theta$, but
by assumption, $\F_1[D]$ acts faithfully on $A$, a contradiction. Thus
Claim 1 follows.
\par
Now suppose $\LL$ is a Lie ideal of $A[D]$ with
$\LL\supset\F_1$ and $\LL\ne\F_1$. Let $n\ge0$ be minimal so that
there exists $u\in\LL\bs\F_1$ with $lev(u)=n$. Fix such $u$.
\par
{\bf Claim 2}.
$n=0$, i.e., $u\in(A\cap\LL)\bs\F_1$.
\par
Suppose conversely $n>0$.
For any $y\in A[D]$, we can decompose $y$ as $y=y^*+y_0$ such that
$y_0\in A$ and all terms in $y^*$ have degree $\ne0$.
So, we can write $u=u^*+u_0$ and $u^*\ne0$.
Since $A[D]$ acts faithfully on $A$, there exists $a\in A$ with
$u^*(a)\ne0$. Set $y=[u,a]\in\LL$.
Using (8), one can calculate that $y_{0}=u^*(a)\ne0$. Thus $[u,a]\ne0$.
Using (8) again and by induction on $lev(u)$,
we see that $y\in\LL$ has leading level $lev(y)\le n-1$ since
$a\in A$. By the minimal choice of $n$, we
must have $y=[u,a]=y_0\in\F_1\subset A$.
We have $A\ne\F{\sc\!}_1{\ssc\,}a+\F_1$,
otherwise $\F[D]$ does not acts faithfully on $A$.
Choose $b\in A\bs(\F{\ssc\!}_1{\ssc\,}a+\F_1)$. Then by replacing $a$ by
$b$ in the above discussion, we must also have $z_0=[u,b]\in\F_1$.
But then
$$
w=[u,ab]=[u,a]b+a[u,b]=y_0 b+z_0 a\in\LL.
\eqno(14)
$$
Since $y_0\ne0$ and $1,a,b$ are $\F_1$-linearly independent
by the choice of $b$, thus $w\notin\F_1$. But $w\in A$ and so
$lev(w)=0<n$,
this contradicts the minimal choice of $n$. Thus Claim 2  follows.
\par
{\bf Claim 3}. $A\subset\LL$.
\par
Let $u$ be as in Claim 2.
Since $u\notin\F_1$, there exists
$\ptl_i,i\in I$ with $\ptl_i(u)\ne0.$
Then for any $x\in A$, we have
$x\ptl_i(u)=[x\ptl_i,u]\in\LL$, i.e., $A\ptl_i(u)\subset A\cap\LL.$
Obviously, $A\ptl_i(u)$ is an ideal of $(A,\cdot)$, let $\MM\subset
A\cap\LL$ be the maximal ideal of $(A,\cdot)$ containing $A\ptl_i(u)$.
Then $\MM\ne0$, and for any $\ptl\in D,x\in\MM$, $\ptl(x)=[\ptl,x]
\in\LL$, thus $\MM+\ptl(\MM)\subset\LL$. But for any $x,y\in\MM,a\in A$,
we have $a\cdot(x+\ptl(y))=ax-\ptl(a)y+\ptl(ay)\in\MM+\ptl(\MM)$, i.e.,
$\MM+\ptl(\MM)$ is an ideal of $(A,\cdot)$. Thus
$\MM+\ptl(\MM)\subset\MM$ by the maximality of $\MM$, and so
$\MM$ is a $ D$-stable ideal of $(A,\cdot)$. But $A$ is $D$-simple
and $\MM\ne0$, we have $\MM=A$. Therefore $A\subset\LL$.
\par
{\bf Claim 4}. If char${\ssc\,}\F=0$, then $\LL=A[D]$.
\par
For a given $\a\in J$ with $|\a|>0$,
note that ${\rm supp}(\a)$ (cf. (7)$\ssc\,$) is finite.
Using induction on the leading degree,
we can suppose that
we have proved that $x\in\LL$ for all monomial $x\in A[D]$ with
$deg(x)<\a$ and ${\rm supp}(deg(x))\subseteq{\rm supp}(\a)$
(note that $\{\b\in J\,|\,\b<\a$ and ${\rm supp}(\b)\subseteq
{\rm supp}(\a)\}$
is a finite set, therefore we can make such an inductive assumption).
For any element $\g\in J$, we denote the largest index $i$ with $\g_i\ne0$
by $\i(\g)$. If $\g=0$, we set $\i(\g)=-\infty.$
For any $k\in I$, we denote $\d^{(k)}\in J$ to be the element such that
$\d^{(k)}_i=\d_{i,k}$ for all $i\in I$.
Let $j=\i(\a)$, then we can write
$$\a=\b+\g\mbox{ and }
\ptl^{(\a)}=\ptl^{(\g)}\ptl^{(\b)}=\ptl_j^{\a_j}\ptl^{(\b)},
\eqno(15)$$
where $\g=\a_j\d^{(j)}$ such that $\i(\b)<j$ and
 ${\rm supp}(\b)\subset{\rm supp}(\a)$
 and ${\rm supp}(\b)\ne{\rm supp}(\a)$.
Choose $v\in A$ with $\ptl_j(v)\ne0$. Then we\vs{-3pt} have
$$
\matrix{
[x\ptl^{(\g+\d^{(j)})},v\ptl^{(\b)}]
\!\!\!\!&=
[x,v\ptl^{(\b)}]\ptl^{(\g+\d^{(j)})}+x[\ptl^{(\g+\d^{(j)})},v\ptl^{(\b)}]
\vs{4pt}\hfill\cr&
=
v[x,\ptl^{(\b)}]\ptl^{(\g+\d^{(j)})}+x[\ptl_j^{\a_j+1},v]\ptl^{(\b)}
\vs{4pt}\hfill\cr&
\dis
=v[x,\ptl^{(\b)}]\ptl^{(\g+\d^{(j)})}+
x\sum_{s=1}^{\a_j+1}({}^{\a_j+1}_{\ \ s})\ptl_j^s(v)\ptl_j^{\a_j+1-s}
\ptl^{(\b)}.
\hfill\cr}
\eqno(16)$$
Observe that
all terms in the summand of the right-hand side except the term
corresponding to $s=1$ are all in $\LL$ by inductive assumption since
they have degree $\g+(1-s)\d^{(j)}+\b<\a$ and their supports are
$\subseteq{\rm supp}(\a)$.
Also the left-hand side is in $\LL$ since
$\b<\a$ and ${\rm supp}(\b)\subset{\rm supp}(\a)$.
As for the first term of the right-hand side,
its leading degree must be $<\a$ because,
either $[x,\ptl^{(\b)}]=0$ if $\b=0$, or else,
$[x,\ptl^{(\b)}]$ is in the subspace spanned by
$A\ptl^{(\tau)}$ for all $\tau\in J$ with
${\rm supp}(\tau)\subset{\rm supp}(\b)$ and
$\tau\le\b-\d^{(k)}$ for some
$k\ge\i(\b)$,
and for all such $\tau$, one has $\tau+\g+\d^{(j)}\le
\b-\d^{(k)}+\g+\d^{(j)}<\b+\g=\a$ since $k\le\i(\b)<j$.
Thus the first term of the right-hand side is in $\LL$
by inductive assumption. This shows that
$x\ptl_j(u)\ptl^{(\a)}\in\LL$.
Let $\MM\subset\{y\in A\,|\,y\ptl^{(\a)}\in\LL\}$
be the maximal ideal of $(A,\cdot)$ containing the ideal
$A\ptl_j(u)\ptl^{(\a)}$ of $(A,\cdot)$,
then as in the proof of Claim 3, $\MM=A$, i.e.,
$A\ptl^{(\a)}\subset\LL.$ This shows that all $x\in A[D]$
with degree $\a$ is in $\LL$. This completes the proof of
our fourth claim, and thus the theorem follows in this case.
\par
Now suppose char${\ssc\,}\F=p>0$.
Observe that for any $\ptl\in {\rm Der}_{\sF}(A)$,
one has $\ptl^p\in{\rm Der}_{\sF}(A)$.
Thus, for any $\ptl_i$ in (1), $\ptl_i^{p^s},s=1,2,\cdots$ are
all derivations of $A$, and they are
A-linearly independent derivations of $A$
because of the faithful action of $A[D]$ on $A$.
In particular, there are infinitely many A-linearly independent
derivations of $A$ in $\F[D]$.
\par
{\bf Claim 5}. $A\ptl\subset\LL$ for all
$\ptl\in{\rm Der}_{\sF}(A) \cap\F[D] $.
\par
Take $\ptl'\in{\rm Der}_{\sF}(A)\cap\F[D] $
such that $\ptl,\ptl'$ are $A$-linearly independent.
By Claim 3, we have
$$
[x\ptl\ptl',a]-x\ptl(\ptl'(a))=
x(\ptl(a)\ptl'+\ptl'(a)\ptl)\in\LL,\ \forall\,x,a\in A.
\eqno(17)$$
Replacing $x$ by $x\ptl(b)$ we obtain
$x\ptl(b)(\ptl(a)\ptl'+\ptl'(a)\ptl)\in\LL$, interchanging $a$ and $b$
and subtracting the two expression, we conclude that
$x(\ptl(b)\ptl'(a)-\ptl(a)\ptl'(b))\ptl\in\LL$. Since $A[D]$ acts
faithfully on $A$, we can choose $a\in A$ with $\ptl'(a)\ne0$, then
$\ptl'(a)\ptl-\ptl(a)\ptl'\ne0$ since $\ptl,\ptl'$ are $A$-linearly
independent, and so, there exists $b\in A$
with $u=\ptl(b)\ptl'(a)-\ptl(a)\ptl'(b)\ne0$. Thus $A u\ptl\subset\LL$.
Let $\MM\subset A$ be the maximal ideal of $(A,\cdot)$ containing
$A u$, then as in the proof of Claim 3,
we have $\MM=A$. Thus $A\ptl\subset\LL$, and the claim follows.
\par
For any $j\in\Z_+,$ we can write $j=\sum_{s\ge 0}j_sp^s$ with $0\le j_s<p$,
and so for any $\ptl\in D$, $\ptl^j=\prod_{s\ge0}(\ptl^{p^s})^{j_s}$,
where all $\ptl^{p^s}$ are derivations of $A$.
Thus, for any $\ptl^{(\a)}\in B$, we can rewrite $\ptl^{(\a)}$ as
$$
\ptl^{(\a)}=\prod_{s=1}^n d_s^{{\ssc\,}m_s},\
1\le m_s\le p-1.
\eqno(18)$$
where $d_s\in\F[D]\cap{\rm Der}_{\sF}(A),1\le s\le n$ are
some $A$-linearly independent derivations.
\par
{\bf Claim 6}. $A\ptl^{(\a)}\subset\LL.$
\par
We shall prove this  by induction on $N_{\a}=\sum_{s=1}^n m_s$
(there may be more than one way to write $\ptl^{(\a)}$ as in (18),
and so $N_{\a}$ may not be uniquely defined. But from the following proof one can
see that this will not affect our inductive step).
By Claim 3 and Claim 5, we have the result in Claim 6 if $N_{\a}\le1$.
So suppose $N_{\a}\ge2$. By inductive assumption,
we can suppose
$$
A\ptl^{(\b)}\subset\LL
\mbox{ for all }\ptl^{(\b)}
\mbox{ with }N_{\b}<N_{\a}.
\eqno(19)$$
Choose $d_{n+1}\in\F[D]\cap{\rm Der}_{\sF}(A)$ such that
$d_1,\cdots,d_{n+1}$ are A-linearly independent. Set $m_{n+1}=1$.
As in (17), for any $a,x\in A$, we have  $[x\ptl^{(\a)}d_{n+1},a] \in\LL$ and
$$
\matrix{
[xd_{n+1} \ptl^{(\a)},a]\!\!\!\!&=
x[d_{n+1},a]\ptl^{(\a)}+xd_{n+1}[\ptl^{(\a)},a]
\vs{4pt}\hfill\cr&\equiv
\dis
xd_{n+1}(a)\ptl^{(\a)}+xd_{n+1}
\sum_{s=1}^nm_sd_s(a)\prod_{r\ne s,1\le r\le n}d_r^{{\ssc\,}m_r}d_s^{{\ssc\,}m_s-1}
\,\,\,( {\rm mod}\,\, {\ssc\,}\LL)
\hfill\cr&\equiv
\dis
x\sum_{s=1}^{n+1}m_sd_s(a)\prod_{r\ne s,1\le r\le n+1}d_r^{{\ssc\,}m_r}d_s^{{\ssc\,}m_s-1}
 \,\,( {\rm mod}\,\, {\ssc\,}\LL),
\,\ \forall\,\,x,a\in A,
\hfill\cr}
\eqno(20)$$
where the second, third equalities follow from the inductive assumption
(19).
Since $d_1,\cdots,$ $d_{n+1}$ are A-linearly independent and $A[D]$
acts faithfully on $A$, using induction on $n$, there exist
$a_1,\cdots,a_{n+1}\in A$ such that
$$
{\rm det}(d_s(a_r))=
\left|\matrix{
d_1(a_1)\!\!&\cdots\!\!&d_1(a_{n+1})\cr
d_2(a_1)\!\!&\cdots\!\!&d_2(a_{n+1})\cr
\cdots\!\!&\cdots\!\!&\cdots\cr
d_{n+1}(a_1)\!\!&\cdots\!\!&d_{n+1}(a_{n+1})\cr
}\right|\in A\bs\{0\}.
\eqno(21)$$
In (20), if we take $a$ to be $a_1,\cdots,a_{n+1}$, and denote
$b_s=[x\ptl^{(\a)}d_{n+1},a_s]\in\LL$ and denote
$y_s=\prod_{r\ne s,1\le r\le n+1}d_r^{{\ssc\,}m_r}d_s^{{\ssc\,}m_s-1}$
for $s=1,\cdots,n+1$, then we obtain equations on $y_s$:
$$
\sum_{s=1}^{n+1}m_sxd_s(a)y_s\equiv b_s\in\LL,\, s=1,\cdots,n+1.
\eqno(22)$$
Note that by (18), $1\le m_s\le p-1$, thus by (21),
we can solve that there exist some $u_s\in A\bs\{0\}$ such
that $x u_s y_s\in\LL$, $1\le s\le n+1$.
In particular, by taking $s=n+1$, since $y_{n+1}=\ptl^{(\a)}$,
we have
$x u_{n+1}\ptl^{(\a)}\in\LL$ for all $x\in A$, i.e., $A u_{n+1}\ptl^{(\a)}
\in\LL$. Hence, as in the proof of Claim 5, we have $A\ptl^{(\a)}
\subset\LL$. This proves the claim, and thus Theorem 1.1 follows.
\qed\par\
\par
{\bf Theorem 1.2} The associative
algebra ${A}[D]$ is simple if and only if
$A$ is $D$-simple and $\F_1[D]$ acts faithfully on $A$.
\par
{\bf Proof.}
``$\Rar$'': Since $\theta$ is a homomorphism of associative algebras,
 $\kn\theta$ is an ideal of $A[D]$. Since $A [D] $ is simple,
if  $\kn\theta=A [D] $, then $A [D] (A)=0$, and in particular
$D(A)=0$, a contradiction with that $D$ is a nonzero subspace
of ${\rm Der}_{\sF}(A)$.
Thus $\kn\theta=0$, i.e.,
$A[D] $ acts faithfully on $A$,  and so $\F_1[D]$ acts
faithfully on $A$.
\par
If $\II$ is a $D$-stable ideal of $A$, then
clearly $\II[D]$ is a ideal of $A[D]$.
If $\II\ne0$, then $\II[D]= A[D]$, and
thus $\II$ must be $A$.
So $A$ is $D$-simple.
\par
``$\Lar$'': Suppose $\II$ is a nonzero ideal of the associative algebra $A[D]$.
Then $\II$ is an ideal of the Lie algebra $A[D]$.
{}From Theorem 1.1, we know that the Lie algebra $\ol{A}[D]$ is simple.
Then $\II=A[D]$ or $\II\subset\F_1$. If $\II\subset\F_1$,
since $\F_1$ is a field,
we must have $\II=A[D]$. Thus, in any case, $\II=A[D]$. So
$A[D]$ is a simple
associative algebra.\qed
\par\
\vs{7pt}\par
{\bf 2 Applications}
\vs{10pt}\par
Using Theorems 1.1, 1.2, one can obtain some interesting
examples of simple (associative or Lie) algebras of Weyl type.
Some are well-known examples, some seem to be new.
\par
{\bf Corollary 2.1}
Let $A=\C[t]$ or $\C[t^{\pm}]$ be the polynomial ring or
Laurent polynomial ring,
$D={\rm span\ssc\,}\{{\ptl\over\ptl t}\}$, then
$
\ol{A}[D]= \C[t,{\ptl\over\ptl t}]/\C\mbox{    and   }
\C[t^{\pm1},{\ptl\over\ptl t}]/\C,
$
are the
well-known Lie algebras of differential operators. More generally,
$$
\C[t_1,\cdots,t_n,
{\ptl\over\ptl t_1},\cdots,{\ptl\over\ptl t_n}]/\C\mbox{ or }
\C[t^{\pm1}_1,\cdots,t^{\pm1}_n,
{\ptl\over\ptl t_1},\cdots,{\ptl\over\ptl t_n}]/\C,
\eqno(23)$$
are simple Lie algebras.
In the first case of (23),
$D$ is a subspace of locally nilpotent derivations,
while in the second case, $\ptl\over\ptl t_i$ can be replaced by
$t_i{\ptl\over\ptl t_i}$ so that $D$ is a subspace of semi-simple
derivations.
The associative algebra $\C[t_1,\cdots,t_n,
{\ptl\over\ptl t_1},\cdots,{\ptl\over\ptl t_n}]$
is the rank $n$ Weyl algebra $A_n$.
\qed\par
{\bf Corollary 2.2}
Let $\ell_1,\ell_2,\ell_3,\ell_4\in\Z_+$ with
$\ell=\ell_1+\ell_2+\ell_3+\ell_4\ge1$.
\vs{-3pt}Let
$$A=\C[t_1,\cdots,t_{\ell_1+\ell_2},t_{\ell_1+\ell_2+1}^{\pm1},
\cdots,t_{\ell_1+\ell_2+\ell_3}^{\pm1},
x_{\ell_1+1}^{\pm1},\cdots,x_{\ell}^{\pm1}],
\vs{-7pt}\eqno(24)$$
and\vs{-7pt} let
$$
D=\{\ptl_i\!=\!{\ptl\over\ptl t_i},
\ptl_j\!=\!{\ptl\over\ptl t_j}\!+\!x_j{\ptl\over\ptl x_j},
\ptl_k\!=\!x_k{\ptl\over\ptl x_k}\,|\,
1\!\le\! i\!\le\!\ell_1\!<\!
j\!\le\!\ell_1\!+\!\ell_2\!+\!\ell_3
\!<\!k\!\le\! \ell\},
\vs{-3pt}\eqno(25)$$
then\vs{-3pt} we obtain nongraded nonlinear simple Lie algebras of Weyl type
$$
\C[t_1,\cdots,t_{\ell_1+\ell_2},
t_{\ell_1+\ell_2+1}^{\pm1},
\cdots,t_{\ell_1+\ell_2+\ell_3}^{\pm1},
x_{\ell_1+1}^{\pm1},\cdots,x_{\ell}^{\pm1},\ptl_1,\cdots,\ptl_\ell]/\C,
\vs{-3pt}\eqno(26)$$
\vs{-3pt}and simple associative algebras of Weyl type
$$
\C[t_1,\cdots,t_{\ell_1+\ell_2},
t_{\ell_1+\ell_2+1}^{\pm1},
\cdots,t_{\ell_1+\ell_2+\ell_3}^{\pm1},
x_{\ell_1+1}^{\pm1},\cdots,x_{\ell}^{\pm1},\ptl_1,\cdots,\ptl_\ell].
\vs{-3pt}\eqno(27)$$
Note that $\ptl_i$ are locally nilpotent, locally finite,
not locally finite or semi-simple
if $1\le i\le\ell_1,$ $\ell_1<i\le\ell_1+\ell_2$,
$\ell_1+\ell_2<i\le\ell_1+\ell_2+\ell_3$ or
$\ell_1+\ell_2+\ell_3<i\le\ell$.
These Lie algebras seem to be new.
\qed\par
{\bf Corollary 2.3} Let $\F$ be a field of arbitrary characteristic.
Let $x_1,x_2,\cdots$ be infinite number of
algebraically independent elements over $\F$.
Let $A=\F(x_1,x_2,\cdots)$ be an extension field of $\F$,
and let $\ptl\in{\rm Der}_{\sF}(\F(x_1,x_2,\cdots))$ such that
$\ptl(x_i)=x_{i+1}$ for $i\ge1$. Set $D=\F\ptl$.
Then $A$ is $D$-simple and $A[D]$ acts faithfully on $A$.
Thus we obtain the simple Lie algebra
$\F(x_1,x_2,\cdots)[\ptl]/\F_1$ and the
simple associative algebra $\F(x_1,x_2,\cdots)[\ptl]$,
where $\F_1=\F$ if char${\ssc\,}\F=0$ or $\F_1=F(x_1^p,x_2^p,\cdots)$
if char${\ssc\,}\F=p>0$.
If char${\ssc\,}\F=2$, this gives an example that
$A[D]/\F_1$ is simple, but the Lie algebra $AD$ considered
in [9] is not simple.
\qed\par
{\bf Corollary 2.4} Let $\F$ be a field of characteristic zero.
Let $\G$ be a multiplicative abelian group and let
$A=\F[\G]$ be the group algebra. Since any $\l\in
{\rm Hom}(\G,\F^+)$ (where $\F^+$ is the additive group $\F$)
gives rise a derivation
$\ptl_\l:\sum_{\a\in\G}f_\a\a\mapsto
\sum_{\a\in\G}f_\a\l(\a)\a$, if we take a \F-subspace
$\D$ of ${\rm Hom}(\G,\F^+)$, then we have the
commutative subalgebra $D=\{\ptl_\l\,|\,\l\in\D\}$
of derivations. If $\G^{\D}=
\{\a\in\G\,|\,\l(\a)=0,\,\forall\,\l\in\D\}
=\{1\}$, then $A$ is $D$-simple and
one can prove that $A[D]$ acts faithfully on $A$, and so we
obtain the simple
Lie algebra $\F[\G,\D]/\F$ and the simple associative algebra
$\F[\G,\D]$.
\qed\par
We would like to conclude our paper with some problems.
\par
{\bf Problems.} 1. When $A[D]\cong A'[D']$ as Lie or associative
algebras?
\par
\hskip 2.2cm 2. What is ${\rm Der}_{\sF}(A[D])$?
\par
\hskip 2.2cm 3. What is $H^2(A[D],\F)$?

\vskip 10pt
{\bf References}
 \par
\begin{description}
\item[{[1]}] Kawamoto N.  Generalizations of Witt algebras over a field of characteristic zero.
Hiroshima Math. J., 1985, 16: 417-462

\item[{[2]}] Osborn J M.  New simple infinite-dimensional Lie algebras of
 characteristic 0. J. Alg., 1996, 185: 820-835

\item[{[3]}] Dokovic D Z,  Zhao K.  Derivations, isomorphisms, and
second cohomology of generalized Witt algebras. Trans.
Amer. Math. Soc., 1998, 350(2): 643-664

\item[{[4]}] Dokovic D Z,  Zhao K.  Generalized Cartan type $W$ Lie algebras
in characteristic zero. J. Alg., 1997, 195: 170-210

\item[{[5]}] Dokovic D Z,  Zhao K.  Derivations, isomorphisms, and
second cohomology of generalized Block  algebras.
Alg. Colloq., 1996, 3(3): 245-272

\item[{[6]}] Osborn J M,   Zhao K.  Generalized Poisson bracket and
Lie algebras
of type H in characteristic 0.  Math. Z., 1999, 230: 107-143

\item[{[7]}] Osborn J M,   Zhao K.  Generalized Cartan type K Lie algebras
in characteristic 0,  Comm. Alg., 1997, 25: 3325-3360

\item[{[8]}] Zhao K. 
Isomorphisms between generalized Cartan type W Lie algebras in
characteristic zero.  Canadian J. Math., 1998,  50: 210-224

\item[{[9]}] Passman D P.  Simple Lie algebras of Witt type. J. Alg.,
1998, 206: 682-692

\item[{[10]}] Jordan D A. On the simplicity of Lie algebras of derivations of
commutative algebras. J. Alg., 2000, 228: 580-585

\item[{[11]}] Xu X. New generalized simple Lie algebras of Cartan type over a
field with characteristic 0.  J. Alg.,  2000, 224: 23-58

\item[{[12]}] Su Y, Xu X, Zhang H.  Derivation-simple algebras
and the structures of Lie algebras of Witt type. J. Alg. in press.

\item[{[13]}] Dixmer J.  Enveloping algebras. North Holland,
Amsterdam. 1977

\item[{[14]}] Zhao K.  Automorphisms of algebras of differential operators.
J. of Capital Normal University, 1994, 1: 1-8

\item[{[15]}] Zhao K.  Lie algebra of derivations of algebras of differential
operators. Chinese Science Bulletin, 1993, 38(10): 793-798

\end{description}

\end{document}